\DeclareFontFamily{U}{wncy}{}
\DeclareFontShape{U}{wncy}{m}{n}{%
   <5>wncyr5%
   <6>wncyr6%
   <7>wncyr7%
   <8>wncyr8%
   <9>wncyr9%
   <10>wncyr10%
   <11>wncyr10%
   <12>wncyr6%
   <14>wncyr7%
   <17>wncyr8%
   <20>wncyr10%
   <25>wncyr10}{}
\DeclareMathAlphabet{\cyrille}{U}{wncy}{m}{n}
\def\sh{
\setlength{\unitlength}{.5 pt}
\begin{picture}(40,20)
\put(10,2){\line(1,0){20}}
\put(10,2){\line(0,1){10}}
\put(20,2){\line(0,1){10}}
\put(30,2){\line(0,1){10}}
\end{picture}}

\newcommand{\dlangle}{\langle\!\langle}
\newcommand{\drangle}{\rangle\!\rangle}

\documentclass[a4paper,11pt]{article}

\usepackage{latexsym}
\usepackage{amsmath,amsfonts,amssymb,amsthm}
\newtheorem{thm}{Theorem}[section]
\newtheorem{lem}[thm]{Lemma}
\newtheorem{exple}[thm]{Example}

\newtheorem{cor}[thm]{Corollary}
\newtheorem{prop}[thm]{Proposition}

\newtheorem{rem}[thm]{Remark}

\input epsf

\title{On exponentials of exponential 
generating series}
\author{Roland Bacher}

\date{Roland.Bacher@ujf-grenoble.fr}

\begin{document}
\maketitle

{\sl Abstract: 
After identification of the algebra of exponential generating series 
with the shuffle
algebra of ordinary formal power series, the exponential map
${\exp}_!:X\mathbb K[[X]]\longrightarrow 1+X\mathbb K[[X]]$
for the associated Lie group with multiplication given by the shuffle
product is well-defined over an arbitrary field $\mathbb K$ by a result
going back to Hurwitz.
The main result of this paper states that
${\exp}_!$ (and its reciprocal map ${\log}_!$) induces a
group isomorphism between the subgroup of rational, respectively algebraic 
series of the additive group $X\mathbb K
[[X]]$ and the subgroup of rational, respectively algebraic series
in the group $1+X\mathbb K[[X]]$ endowed with the shuffle product,
if the field $\mathbb K$ is a subfield of
the algebraically closed field $\overline{\mathbb F}_p$ of
characteristic $p$.}

\section{Introduction}

The equality
\begin{eqnarray}\label{formshuffle1}
\left(\sum_{n=0}^\infty \alpha_n\frac{X^n}{n!}\right)
\left(\sum_{n=0}^\infty \beta_n\frac{X^n}{n!}\right)=
\sum_{n=0}^\infty \sum_{m=0}^n{n+m\choose n}
\alpha_n\beta_m\frac{X^{n+m}}{(n+m)!}
\end{eqnarray}
shows that we can define an algebra structure on the vector space
$$\mathcal E(\mathbb K)=\left\{\sum_{n=0}^\infty \alpha_n \frac{X^n}{n!}
\ \vert\ \alpha_0,\alpha_1,\ldots\in\mathbb K\right\}$$
of formal exponential generating series with 
coefficients $\alpha_0,\alpha_1,\dots$
in an arbitrary field or ring $\mathbb K$. For the sake of simplicity 
we work in the sequel only over fields. The expression
$\alpha_n/n!$ should be considered formally since the numerical
value of $n!$ is zero over a field of positive characteristic 
$p\leq n$.

Motivation for this work is given by the fact that formula (\ref{formshuffle1})
allows to define the shuffle product 
$$\sum_{n=0}^\infty \gamma_n X^n=\left(\sum_{n=0}^\infty \alpha_n X^n\right)
\sh\left(\sum_{n=0}^\infty \beta_n X^n\right)$$
of two formal power series $\sum_{n=0}^\infty \alpha_n X^n$ and 
$\sum_{n=0}^\infty \beta_n X^n$ by setting 
\begin{eqnarray}\label{defshuffleprod}
\gamma_n=\sum_{k=0}^n{n\choose k}\alpha_k\beta_{n-k}\ .
\end{eqnarray}
The definition of the shuffle product
arises in the theory of divided powers, see eg. 
\cite[Definition 3.1]{BO}. The main properties needed 
in this paper are however already 
in Hurwitz, see \cite{H}.
I have the impression that the main results of the present
paper, given by Theorem \ref{thmmainrat} and \ref{thmmainalg}
do not fit very well into the theory of divided powers: They are 
based on an
interplay between ordinary power series (used for defining 
rationality and algebraicity) and exponential power series
(used for defining an analogue of the exponential map in 
positive characteristic). A special instance of this exponential 
map is a standard ingredient for divided powers, 
see \cite[Appendix A, Proposition A1]{BO},
but ordinary formal power series do not seem to play a significant 
role there.

Definition \ref{defshuffleprod} 
is also a particular case of a shuffle product defined more generally 
for formal power series in several non-commuting variables.
The associated shuffle-algebras arise 
for example in the study of free Lie algebras \cite{Re}, Hopf algebras and 
polyzetas \cite{Z}, \cite{Ca}, 
formal languages \cite{BR}, etc.

I became interested in this subject through the study of 
the properties of the algebra of recurrence matrices, a subset of 
sequences of matrices displaying a kind of self-similarity structure
used in \cite{B1} and \cite{B} for studying reductions  
of the Pascal triangle modulo suitable Dirichlet characters.
Such recurrence matrices are closely related to automata groups and
complex dynamical systems, see for example \cite{N} for details. 
Over a finite field, they can be identified 
with rational formal power series in several non-commuting variables 
(the underlying algebras are however very different) and it is thus 
natural to investigate properties of other possible products preserving 
these sets. The main results of this paper, given by 
Theorem \ref{thmmainrat} and \ref{thmmainalg} (and their 
effective analogues, Theorem \ref{thmrateff} and \ref{thmalgeff}),
deal with properties of 
the shuffle product for formal power series in one variable which have 
gone unnoticed in the existing literature, as far as I am aware.


We denote by
$$\mathfrak m_{\mathcal E}=
\left\{\sum_{n=1}^\infty \alpha_n \frac{X^n}{n!}\ \vert\
\alpha_1,\alpha_2,\dots\in\mathbb K\right\}\subset \mathcal E(\mathbb K)$$
the maximal ideal of the local algebra $\mathcal E(\mathbb K)$.
A straightforward computation already known to Hurwitz, see \cite{H},
shows that $a^n/n!$ is always well-defined
for $a\in \mathfrak m_{\mathcal E}$. 
Endowing $\mathbb K$ with the discrete topology and
$\mathcal E(\mathbb K)$ with the topology given by 
coefficientwise convergency, the functions
$$\exp(a)=\sum_{n=0}^\infty \frac{a^n}{n!}\hbox{ and }
\log(1+a)=-\sum_{n=1}^\infty \frac{(-a)^n}{n}$$
are always defined for $a\in\mathfrak m_{\mathcal E}$.

Switching back to ordinary generating series 
$$A=\sum_{n=1}^\infty \alpha_nX^n,\ B=\sum_{n=1}^\infty\beta_n X^n 
\in\mathfrak m$$
contained in the maximal ideal $\mathfrak m=X\mathbb K[[X]]$,
of (ordinary) formal power series, we write 
$${\exp}_!(A)=1+B$$
if 
$$\exp\left(\sum_{n=1}^\infty \alpha_n\frac{X^n}{n!}\right)=
1+\sum_{n=1}^\infty \beta_n\frac{X^n}{n!}\ .$$
It is easy to see that ${\exp}_!$ defines a one-to-one
map between $\mathfrak m$ and $1+\mathfrak m$ with reciprocal map
$$1+B\longmapsto A={\log}_!(1+B)\ .$$

It satisfies $${\exp}_!(A+B)={\exp}_!(A)
\sh {\exp}_!(B)$$
for all $A,B\in\mathfrak m$
where the shuffle product
$$\left(\sum_{n=0}^\infty \alpha_nX^n\right)\sh
\left(\sum_{n=0}^\infty \beta_n X^n\right)=
\sum_{n,m=0}^\infty {n+m\choose n}\alpha_n\beta_m X^{n+m}$$
corresponds to the ordinary 
product of the associated exponential generating series. 
The map ${\exp}_!$ defines thus an isomorphism between
the additive group $(\mathfrak m,+)$ and
the {\it special shuffle-group} $(1+\mathfrak m,\sh)$ with group-law
given by the shuffle-product. It coincides with the familiar
exponential map from the Lie algebra  $\mathfrak m$ into the 
special shuffle group, considered as an infinite-dimensional Lie group.

The paper \cite{Fl} of Fliess implies that rational,
respectively algebraic elements form a subgroup in 
$(1+\mathfrak m,\sh)$ if one works over a subfield of 
$\overline{\mathbb F}_p$. It is thus natural to consider the corresponding 
subgroups (under the reciprocal map 
$\log_!$ of the Lie-exponential $\exp_!:
\mathfrak m\longmapsto 1+\mathfrak m$) 
in the isomorphic additive group $(\mathfrak m,+)$ forming 
the Lie algebra of $(1+\mathfrak m,\sh)$. The answer which is the main result
of this paper is surprisingly simple: The corresponding subgroup
is exactly the subgroup of all rational, respectively algebraic
elements in the additive group $\mathfrak m$. We have thus:

\begin{thm} \label{thmmainrat}
Let $\mathbb K$ be a subfield of the algebraically closed
field $\overline{\mathbb F}_p$ of positive characteristic 
$p$. Given a series $A\in\mathfrak m=X\mathbb K[[X]]$ the 
following two assertions are equivalent:

\ \ (i) $A$ is rational.

\ \ (ii) ${\exp}_!(A)$ is rational.
\end{thm}

Theorem \ref{thmmainrat} fails in
characteristic zero: The series 
$${\log}_!(1-X)=-\sum_{n=1}^\infty (n-1)!X^n$$ 
is obviously transcendental. (This series shows also that 
Theorem \ref{thmmainalg} fails to hold in characteristic zero.)

\begin{exple} The Bell numbers $B_0,B_1,B_2,\dots$, see pages 45,46 in 
\cite{Comt} or Example 5.2.4 in 
\cite{St2}, are the natural integers
defined by 
$$\sum_{n=0}^\infty B_n\frac{x^n}{n!}=e^{e^x-1}$$
and have combinatorial interpretations.

Since $e^x-1$ is the exponential generating series of the sequence
$0,1,1,\dots$, we have $\sum_{n=0}^\infty B_n x^n=\exp_!(x/(1-x))$
for the ordinary generating series
$$\sum_{n=0}^\infty B_nx^n=1+x+2x^2+5x^3+15x^4+52x^5+203x^6+877x^7+
4140x^8+\dots$$
of the Bell numbers.

The reduction of $\sum_{n=0}^\infty B_n x^n$ modulo a prime $p$ is thus
always a rational element of $\mathbb F_p[[x]]$.
A few such reductions are
$$\frac{1}{1+x+x^2}\pmod 2,\ \frac{1+x+x^2}{1-x^2-x^3}\pmod 3,\ 
\frac{1+x+2x^2-x^4}{1-x^4-x^5}\pmod 5\ .$$
\end{exple}

\begin{thm} \label{thmmainalg}
Let $\mathbb K$ be a subfield of the algebraically closed
field $\overline{\mathbb F}_p$ of positive characteristic 
$p$. Given a series $A\in\mathfrak m=X\mathbb K[[X]]$ the 
following two assertions are equivalent:

\ \ (i) $A$ is algebraic.

\ \ (ii) ${\exp}_!(A)$ is algebraic.
\end{thm}
 
Theorem \ref{thmmainrat} and \ref{thmmainalg} are the main results of this
paper and can be restated as follows.

\begin{cor} Over a subfield $\mathbb K\subset \overline{\mathbb F}_p$,
the group isomorphism 
$${\exp}_!:(\mathfrak m,+)\longrightarrow 
(1+\mathfrak m,\sh)$$
restricts to an ismorphism between the subgroups of rational 
elements in $(\mathfrak m,+)$ and in
$(1+\mathfrak m,\sh)$. 

It restricts also to an isomorphism between the subgroups of 
algebraic elements in $(\mathfrak m,+)$ and in $(1+\mathfrak m,\sh)$.

In particular, the subgroup of rational, respectively algebraic elements
in the shuffle group $(1+\mathfrak m,\sh)$ is a Lie-group whose Lie algebra
(over $\mathbb K\subset \overline{\mathbb F}_p$) is given by the 
additive subgroup of all rational, respectively algebraic elements in
$(\mathfrak m,+)$.
\end{cor}

Theorem \ref{thmmainrat} and \ref{thmmainalg} 
can be made more precise as follows.

Given a rational series $A\in\mathbb K[[X]]$
represented by a reduced fraction $f/g$
where $f,g$ with $g\not=0$ are two coprime polynomials of
degree $\deg(f)$ and $\deg(g)$, we set 
$\parallel A\parallel=\max(1+\deg(f),\deg(g))$, see also Proposition 
\ref{propcharrat} for a well-known equivalent description of 
$\parallel A\parallel$.

\begin{thm} \label{thmrateff} We have
$$\parallel {\exp}_!(A)\parallel \leq p^{q^{\parallel A\parallel}} 
\hbox{ and }\parallel {\log}_!(1+A)\parallel \leq 1+
\parallel 1+A\parallel^p$$
for a rational series $A$ in $\mathfrak m\subset 
\overline{\mathbb F}_p[[X]]$ having all its coefficients in a
finite subfield $\mathbb F_q\subset\overline{\mathbb F}_p$ containing 
$q=p^e$ elements.
\end{thm}

The bounds for $\log_!$ (and the analogous bounds in the 
algebraic case) can be improved,
see Proposition \ref{propsra}. 

Theorem \ref{thmrateff} could be called an effective version of Theorem
\ref{thmmainrat}: Given a rational series represented by $f/g\in
\mathfrak m$ with $f,g\in\overline{\mathbb F}_p[X]$, Theorem 
\ref{thmmainrat} ensures the existence of polynomials $u,v$
such that $\exp_!(f/g)=u/v$. Theorem \ref{thmrateff}
shows that $u$ and $v$ are of degre at most $p^{q^{\parallel f/g\parallel}}$.
They can thus be recovered as suitable Pad\'e approximants from the
series developement of $\exp_!(f/g)$ up to order 
$2p^{q^{\parallel f/g\parallel}}$. Experimentally, the number $\parallel
\exp_!(A)\parallel$ is generally much smaller.

Since the bounds for $\log_!$ are better than for 
$\exp_!$, the determination of
the rational series $B=\exp_!(A)$ with $A\in\mathfrak m$ 
rational is best done as follows: Start by ``guessing'' the rational series 
$B$ and check (or improve the guess for $B$ in case of failure)
that $A=\log_!(B)$ using the bounds for 
$\log_!$. 

Given a prime $p$ and a formal power series 
$C=\sum_{n=0}^\infty \gamma_nX^n$ in
$\mathbb K[[X]]$ with coefficients in a subfield $\mathbb K$ of
$\overline{\mathbb F}_p$, we define
for $f\in\mathbb N,\ k\in\mathbb N,\ k<p^f$ the series
$$C_{k,f}=\sum_{n=0}^\infty 
\gamma_{k+np^f}X^n\ .$$
The vector space $\mathcal K(C)=\mathbb K C+\sum_{k,f}\mathbb KC_{k,f}$
spanned by $C$ and by all series of the form $C_{k,f},\
k\in\{0,\dots,p^f-1\},f\in \{1,2,\dots\}$ is called the {\it $p-$kernel}
of $C$. We denote its dimension by $\kappa(C)=\dim(\mathcal K(C))$.

Algebraic series in $\mathbb K[[X]]$ for $\mathbb K$ a subfield
of $\overline{\mathbb F}_p$ are characterised by a
Theorem of Christol (see Theorem 12.2.5 in \cite{AS}) stating that
a series $C$ in $\overline{\mathbb F}_p[[X]]$ is algebraic
if and only if its $p-$kernel $\mathcal K(C)$ is of finite
dimension $\kappa(C)<\infty$.
We have $\kappa(A+B) \leq \kappa(A)+
\kappa(B)$ and an algebraic series $A\in\overline{\mathbb F}_p[[X]]$
has a minimal polynomial of degree at most $p^{\kappa(A)}$
with respect to $A$.

\begin{thm} \label{thmalgeff} We have
$$\kappa({\exp}_!(A)) \leq q^{\kappa(A)-1}p^{q^{\kappa(A)}}
\hbox{ and }\kappa(
{\log}_!(1+A))
\leq 1+4(\kappa(1+A))^p$$
for a non-zero algebraic series $A$ in $\mathfrak m\subset 
\overline{\mathbb F}_p[[X]]$ having all its coefficients in a
finite subfield $\mathbb F_q\subset\overline{\mathbb F}_p$ containing 
$q=p^e$ elements.
\end{thm}

Considerations similar to those made after Theorem \ref{thmrateff} 
are valid and Theorem \ref{thmalgeff} can be turned into an algorithmically 
effective version of Theorem \ref{thmmainalg}.

A map $\mu:\mathcal V\longrightarrow \mathcal W$ between 
two $\mathbb K-$vector spaces is a homogeneous 
form of degree $d$ if $l\circ \mu:\mathcal V\longrightarrow \mathbb K$
is homogeneous of degree $d$ (given by a homogeneous polynomial of degre $d$
with respect to coordinates) for every linear form 
$l:\mathcal W\longrightarrow \mathbb K$.

A useful ingredient for proving Theorems \ref{thmmainrat}, 
\ref{thmmainalg} and their effective versions is the following 
characterisation of ${\log}_!$:

\begin{prop} Over a field $\mathbb K\subset \overline{\mathbb F}_p$,
the application ${\log}_!:1+\mathfrak m\longrightarrow 
\mathfrak m$ extends to a homogeneous form of degree $p$
from $\mathbb K[[X]]$ into $\mathfrak m$.
\end{prop}

\begin{exple} In characteristic $2$, we have
$${\log}_!\left(\sum_{n=0}^\infty\alpha_nX^n\right)=
\sum_{n=0}^\infty \alpha_{2^n}^2X^{2^{n+1}}+\sum_{0\leq i<j}
{i+j\choose i}\alpha_i\alpha_jX^{i+j}$$
for $\sum_{n=0}^\infty \alpha_nX^n$ in $1+X{\overline{\mathbb F}}_2[[X]]$.
\end{exple}

\begin{rem} Defining $f_!$ as $$f_!\left(\sum_{n=1}^\infty 
\alpha_nX^n\right)=\sum_{n=1}^\infty \beta_nX^n$$
if 
$$f\left(\sum_{n=1}^\infty \alpha_n\frac{X^n}{n!}\right)=
\sum_{n=1}^\infty \beta_n\frac{X^n}{n!}$$
Theorems \ref{thmmainrat}, \ref{thmmainalg}, \ref{thmrateff}
and \ref{thmalgeff} have analogs for the functions
$\mathop{sin}_!$ and $\mathop{tan}_!$ (and for their reciprocal
functions $\mathop{arcsin}_!$ and $\mathop{arctan}_!$).
\end{rem}

The rest of the paper is organised as follows. 
In Sections \ref{sectrat}-\ref{sectlog} we recall
a few definitions and well-known facts which are essentially 
standard knowledge in the theory of divided powers, 
see \cite{BO} or the original work \cite{R1,R2}.
Section \ref{sectproofs} contains the proofs for all
results mentionned above. 

In a second part, starting at Section \ref{sectionpows},
we generalise Theorems \ref{thmmainrat} and
\ref{thmrateff} to formal power series in several non-commuting variables.

\section{Rational and algebraic elements in $\mathbb K[[X]]$}
\label{sectrat}

This section recalls a few well-known facts concerning rational
and algebraic elements in the algebra $\mathbb K[[X]]$ of
formal power series.

We denote by $\tau:\mathbb K[[X]\longrightarrow \mathbb K[[X]]$
the shift operator
$$\tau\left(\sum_{n=0}^\infty \alpha_nX^n\right)=\sum_{n=0}^\infty
\alpha_{n+1}X^n$$
acting on formal power series. The following well-known result
characterises rational series:

\begin{prop} \label{propcharrat}
A formal power series $A=\sum_{n=0}^\infty \alpha_nX^n$ of 
$\mathbb K[[X]]$ is rational
if and only if the series $A,\tau(A),\tau^2(A),\dots,\tau^k(A)=
\sum_{n=0}^\infty \alpha_{n+k}X^n,\dots$
span a finite-dimensional vector space in $\mathbb K[[X]]$.

More precisely, the vector space spanned by 
$A,\tau(A),\tau^2(A),\dots,\tau^i(A),\dots$
has dimension $\parallel A\parallel=\max(1+\deg(f),\deg(g))$
if $f/g$ with $f,g\in\mathbb K[X]$ is a reduced expression
of a rational series $A$.
\end{prop}

The function $A\longmapsto\parallel A\parallel$ satisfies
the inequality 
$$\parallel A+B\parallel \leq \parallel A\parallel+\parallel B\parallel$$ 
for rational series $A,B$ in $\mathbb K[[X]]$. As a particular case we have
$$\parallel A\parallel-1\leq \parallel 1+A\parallel \leq \parallel A
\parallel +1\ .$$

Given a prime $p$ and a formal power series 
$C=\sum_{n=0}^\infty \gamma_nX^n$ in
$\overline{\mathbb F}_p[[X]]$ we denote by $\kappa(C)\in\mathbb N
\cup\{\infty\}$ the dimension of its $p-$kernel
$$\mathcal K(C)=\mathbb KC+\sum_{f,k}\overline{\mathbb F}_pC_{k,f}$$
spanned $C$ and by all series of the form 
$$C_{k,f}=\sum_{n=0}^\infty \gamma_{k+np^f}X^n$$
with $k\in\mathbb N$ such that $k<p^f$ for $f\in\{1,2,\dots\}$.

Algebraic series of $\mathbb K[[X]]$ for $\mathbb K$ a subfield
of the algebraic closure $\overline{\mathbb F}_p$ of finite
prime characteristic $p$ are characterised by the following
Theorem of Christol (see \cite{Chr} or Theorem 12.2.5 in \cite{AS}):

\begin{thm} \label{thmcharalg} 
A formal power series $C=\sum_{n=0}^\infty \gamma_nX^n$ of
$\overline{\mathbb F}_p[[X]]$ is algebraic if and only if the dimension
$\kappa(C)=\dim(\mathcal K(C))$ of its $p-$kernel $\mathcal K(C)$ is finite.
\end{thm}

Finiteness of $\kappa(C)$ amounts to recognisability of $C$ 
which has the following well-known consequence.

\begin{cor} \label{coralgfin}
An algebraic series of $\overline{\mathbb F}_p[[X]]$
has all its coefficients in a finite subfield of 
$\overline{\mathbb F}_p$.
\end{cor}

\begin{prop} \label{propKtau}
Let $C=\sum_{n=0}^\infty \gamma_nX^n$ be an algebraic series 
with coefficients in a subfield $\mathbb K\subset \overline{\mathbb F}_p$.

\ \ (i) We have 
$$\mathcal K(\tau(C))\subset\mathcal K(C)
+\tau(\mathcal K(C))$$
which implies 
$$\kappa(\tau(C))\leq 2\kappa(C)\ .$$ 

\ \ (ii) We have 
$$\mathcal K(C)\subset \mathbb K+\mathcal K(\tau(C))+
X\mathcal K(\tau(C))$$
which implies 
$$\kappa(C)\leq 1+2\kappa( \tau(C))\ .$$
\end{prop}

{\bf Proof} Assertion (i) 
follows from an iterated application of the easy computations
$$(\tau(C))_{k,1}=C_{k+1,1}$$
if $0\leq k< p-1$ and 
$$(\tau(C))_{p-1,1}=\tau(C_{0,1})\ .$$

The proof of assertion (ii) is similar.\hfill$\Box$

\section{The shuffle algebra}\label{sectshufflpr}

This section recalls mostly well-known results concerning shuffle
products of elements in the set $\mathbb K[[X]]$ of formal power
series over a commutative field $\mathbb K$ which is arbitrary
unless specified otherwise. 

The {\it shuffle product} 
$$A\sh B=C=\sum_{n=0}^\infty \gamma_nX^n$$
of $A=\sum_{n=0}^\infty \alpha_nX^n$ and $B=\sum_{n=0}^\infty \beta_nX^n$
is defined by
$$\gamma_n=\sum_{k=0}^n {n\choose k}\alpha_k\beta_{n-k}$$
and corresponds to the usual product  
$ab=c$ of the associated exponential generating series
$$a=\sum_{n=0}^\infty \alpha_n\frac{X^n}{n!},\ 
b=\sum_{n=0}^\infty \beta_n\frac{X^n}{n!},\ 
c=\sum_{n=0}^\infty \gamma_n\frac{X^n}{n!}\ .$$
The {\it shuffle algebra} is the algebra $(\mathbb K[[X]],\sh)$
obtained by endowing the vector space $\mathbb K[[X]]$ of ordinary 
generating series with the shuffle product. By construction,
the shuffle algebra
is isomorphic to the algebra $\mathcal E(\mathbb K)$
of exponential generating series.
In characteristic zero, the 
trivial identity 
$$\sum_{n=0}^\infty \alpha_nX^n=\sum_{n=0}^\infty (n!\alpha_n)\frac{X^n}{n!}$$
gives an isomorphism between the usual algebra $\mathbb K[[X]]$
of ordinary generating series and the shuffle algebra $(\mathbb K[[X]],\sh)$.

The identity $\left(\sum_{n\geq 0}\lambda^nX^n\right)
\sh\left(\sum_{n\geq 0}\mu^nX^n\right)=\sum_{n\geq 0}
(\lambda+\mu)^nX^n$, equivalent to
$e^{\lambda X}e^{\mu X}=e^{(\lambda+\mu)X}$ implies that
the convergency radius of the 
shuffle product of two complex series with strictly positive 
convergency radii $\rho_1,\rho_2$ is at least the harmonic mean
$1/\left(\frac{1}{\rho_1}+\frac{1}{\rho_2}\right)$ of 
$\rho_1$ and $\rho_2$.

\begin{prop} \label{propshiftderiv}
The shift operator $\tau(\sum_{n=0}^\infty \alpha_nX^n)=
\sum_{n=0}^\infty \alpha_{n+1}X^n$ acts as a derivation on the 
shuffle algebra.
\end{prop}

{\bf Proof} The map $\tau$ is clearly linear. The
computation
$$\begin{array}{l}
\displaystyle \tau\left(\sum_{i,j\geq 0}
{i+j\choose i}\alpha_i\beta_jX^{i+j}\right)=
\sum_{i,j\geq 0}{i+j\choose i}\alpha_i\beta_jX^{i+j-1}=\\
\displaystyle \quad =
\sum_{i,j\geq 0}
\left({i+j-1\choose i-1}+{i+j-1\choose j-1}\right)\alpha_i\beta_j
X^{i+j-1}\end{array}$$
shows that $\tau$ satisfies the Leibniz rule
$\tau(A\sh B)=\tau(A)\sh B+A\sh \tau(B)$.\hfill$\Box$ 

Proposition \ref{propshiftderiv} is trivial and well-known
in characteristic zero: 
the usual derivation $d/dX$ acts obviously as the shift operator
on the algebra $\mathcal E(\mathbb K)$ of exponential generating series
over a field of characteristic zero.

The following two results seem to be due to Fliess, cf. Proposition
6 in \cite{Fl}.

\begin{prop} \label{proprat} Shuffle products of rational power series are
rational. 

More precisely, we have 
$$\parallel A\sh B\parallel \leq \parallel A\parallel \ 
\parallel B\parallel$$
for two rational series $A,B$ in $\mathbb K[[X]]$.
\end{prop}

{\bf Proof} Proposition \ref{propshiftderiv} implies 
$\tau^n\left(A\sh B\right)=\sum_{k=0}^n{n\choose k}
\tau^k(A)\sh \tau^{n-k}(B)$.
The series $\tau^n\left(A\sh B\right)$ belongs thus to the vector space
spanned by shuffle products with factors in the vector spaces
$\sum_{n\geq 0}\mathbb K\tau^n(A)$ and $\sum_{n\geq 0}\mathbb K\tau^n(B)$.
This implies the inequality. 
Proposition \ref{propcharrat} ends the proof.
\hfill$\Box$ 

\begin{prop} \label{propalg} Shuffle products of algebraic series in 
$\overline{\mathbb F}_p[[X]]$ are algebraic.

More precisely, we have 
$$\kappa(A\sh B)\leq \kappa(A) \
\kappa(B)\ .$$
\end{prop}

{\bf Proof} Denoting as in Section \ref{sectrat} by $C_{k,f}$
the series 
$$C_{k,f}=\sum_{n=0}^\infty \gamma_{k+np^f}X^n$$
associated to a series $C=\sum_{n=0}^\infty \gamma_n X^n$
and by $\kappa(C)$ the dimension of the vector space
$\mathcal K(C)=\mathbb KC
+\sum_{k,f}\overline{\mathbb F}_pC_{k,f}$, Lucas's identity
(see \cite{Lu})
$${n\choose k}\equiv \prod_{i\geq 0}{\nu_i\choose \kappa_i}\pmod p$$
for $n=\sum_{i\geq 0}\nu_ip^i,\ k=\sum_{i\geq 0}\kappa_ip^i$ with 
$\nu_i,\kappa_i\in\{0,\dots,p-1\}$ implies
$$\left(A\sh B\right)_{k,1}=\sum_{i=0}^k{k\choose i}
A_{i,1}\sh B_{k-i,1}$$
for $k=0,\dots,p-1$. Iteration of this formula shows that 
$\left(A\sh B\right)_{k,f}$ (for arbitrary $k,f\in\mathbb N$
such that $k<p^f$) belongs to the vector space spanned by shuffle 
products with factors in the vector spaces
$\mathcal K(A)$ and $\mathcal K(B)$ of dimension
$\kappa( A)$ and $\kappa(B)$.

Christol's Theorem (Theorem \ref{thmcharalg}) ends the proof.
\hfill$\Box$

\begin{rem} \label{remrefinement} {\rm Given a subfield $\mathbb K$ of 
$\overline{\mathbb F}_p$
let $\mathcal A\subset \mathbb K[[X]]$ denote a 
vector space of finite dimension $a=\dim(\mathcal A)$ 
containing the $p-$kernel $\mathcal K(A)$ of every element 
$A\in\mathcal A$.

We consider an element $B=A_1\sh A_2\sh\cdots \sh A_k$ given 
by the shuffle product of $k$ series $A_1,\dots,A_k\in\mathcal A$.
Expressing all elements $A_1,A_2,\dots$ as linear combinations of 
elements in a fixed basis of $\mathcal A$ and using commutativity
of the shuffle product, the proof of Proposition \ref{propalg} shows
that the inequality $\kappa(B)\leq \kappa(A_1)\kappa(A_2)\cdots\leq
a^k=(\dim(\mathcal A))^k$ can be improved to 
$$\kappa(B)\leq {k+a-1\choose a-1}$$
where the binomial coefficient 
${k+a-1\choose a-1}$ encodes the dimension of the vector space of 
homogeneous polynomials of degree $k$ in $a$ (commuting) variables
$X_1,X_2,\dots,X_a$.}
\end{rem}

\section{The special shuffle-group}

We call the group of units of the shuffle algebra $(\mathbb K[[X]],\sh)$
the {\it shuffle-group}. Its elements are given by the set $\mathbb K^*
+X\mathbb K[[X]]$ underlying the multiplicative unit group.
The shuffle-group is the direct product of the unit group $\mathbb K^*$
of $\mathbb K$ with the {\it special shuffle-group} $(1+X\mathbb K[[X]],
\sh)$.

The inverse in the shuffle group of 
$1-A\in (1+X\mathbb K[[X]],\sh)$ is given by
$$\sum_{n=0}^\infty A^{\sh^n}=1+A+A\sh A+A\sh A\sh
A+\dots$$
where $A^{\sh^0}=1$ and $A^{\sh^{n+1}}=A\sh A^{\sh^n}$ for $n\geq 1$.

The trivial identity $X\sh X^n={n+1\choose 1}X^{n+1}=(n+1)X^{n+1}
\in \mathbb K[[X]]$ implies
$(1-X)\sh\left(\sum_{n=0}^\infty n!X^n\right)=1$. 
Invertible rational (analytical) power series have thus generally a 
transcendental (non-analytical) shuffle-inverse over the complex numbers.

\begin{prop} \label{propshgrFp} 
The special shuffle-group $(1+X\mathbb K[[X]],\sh)$ is isomorphic to an
infinite-dimensional $\mathbb F_p-$vector space if 
the field $\mathbb K$ is of positive characteristic $p$.
\end{prop}

Proposition \ref{propshgrFp} shows that
$(1+X\mathbb K[[X]],\sh)$ is
not isomorphic to the multiplicative group structure on $1+X\mathbb
K[[X]]$ if $\mathbb K$ is of positive characteristic.

{\bf Proof of Proposition \ref{propshgrFp}} Follows from
the fact that ${\exp}_!$ is a group isomorphism between
the $\mathbb F_p-$vector space $\mathfrak m$ and the special
shuffle group. 
\hfill$\Box$

Proposition \ref{propshgrFp} follows also as a special case from 
Proposition \ref{propshgrFpgen}. This yields a different proof 
which is not based on properties of $\exp_!$.

\begin{rem} \label{remshfflegrpnoalgelts}
One can show that a rational fraction 
$A\in 1+X\mathbb C[[X]]$
has a rational inverse for the shuffle-product if and only if
$A=\frac{1}{1-\lambda X}$ with
$\lambda\in\mathbb C$. (Compute $A\sh B=1$
using the decomposition into simple fractions of the rational series 
$A,B$.)
\end{rem}

\section{The exponential and the logarithm for 
exponential generating functions}

Hurwitz showed that $\frac{1}{k!}a^k$ is well-defined for
$a\in\mathfrak m_{\mathcal E}$ with coefficients in an arbitrary
field or commutative ring, see Satz 1 in \cite{H}. We give a
different proof of this fact which implies that $\exp_!$
and $\log_!$ are well-defined over fields of positive characteristic.

\begin{prop} \label{propcoeffintegr} For all natural numbers $j,k\geq 1$,
the set $\{1,\dots,jk\}$ can be partitioned in 
exactly $$\frac{(jk)!}{(j!)^k\ k!}$$
different ways into $k$ unordered disjoint subsets of $j$
elements.

In particular, the rational number $(jk)!/((j!)^k\ k!)$
is an integer for all natural numbers $j,k$ such that $j\geq 1$.
\end{prop}

{\bf Proof} The multinomial coefficient $(jk)!/(j!)^k$ counts
the number of ways of partitioning $\{1,\dots,jk\}$ into
an ordered sequence of $k$ disjoint subset containing all $j$ elements. 
Dividing by $k!$ removes the order on these $k$ subsets.

This proves that the formula defines an integer for all $j,k\geq 1$ and
integrality holds also obviously for $k=0$ and $j\geq 1$.
\hfill$\Box$

\begin{rem} {\rm
A slightly different proof of Proposition \ref{propcoeffintegr} follows from
the observation that $(jk)!/((j!)^k\ k!)$ is the index in the symmetric group
over $jk$ elements of the subgroup formed by all permutations stabilising
a partition of the set $\{1,\dots,jk\}$ into $k$ disjoint subsets of $j$
elements.

A different proof is given by the formula $\frac{(jk)!}{(j!)^kk!}=\prod_{n=1}^k{nj-1\choose j-1}$, easily shown using induction on $k$, see \cite[Section 3]{BO} (which contains a small misprint).}
\end{rem}

\begin{prop} For any natural integer $k\in\mathbb N$, there exists
polynomials $P_{k,n}\in\mathbb N[\alpha_1,\dots,\alpha_n]$ such that 
$$\frac{1}{k!}\left(\sum_{n=1}^\infty \alpha_n\frac{X^n}{n!}\right)^k=
\sum_{n=0}^\infty P_{k,n}(\alpha_1,\alpha_2,\dots,\alpha_n)\frac{X^n}{n!}
\ .$$
\end{prop}

{\bf Proof} The contribution of a monomial 
$$\alpha_1^{j_1}\alpha_2^{j_2}
\dots\alpha_s^{j_s}\frac{X^{\sum_{i=1}^sij_i}}{\left(\sum_{i=1}^sij_i
\right)!}$$
with $j_1+j_2+\dots+j_s=k$ 
to $(1/k!)\left(\sum_{n=1}^\infty \alpha_nX^n/n!\right)^k$
is given by
$$\frac{1}{k!}\frac{k!}
{(j_1)!(j_2)!\cdots (j_s)!}\frac{\left(\sum_{i=1}^s ij_i\right)!}{
\prod_{i=1}^s (i!)^{j_i}}$$
$$=\left(\prod_{i=1}^s\frac{(ij_i)!}{(i!)^{j_i}(j_i)!}\right)
\frac{\left(\sum_{i=1}^sij_i\right)!}{\prod_{i=1}^s(i j_i)!}$$
and the last expression is a product of a natural integer by Proposition 
\ref{propcoeffintegr} and of a multinomial coefficient. 
It is thus a natural integer.
\hfill$\Box$

\begin{cor} For $a=\sum_{n=1}^\infty \alpha_n \frac{X^n}{n!}$ the formulae
$$\exp\left(\sum_{n=1}^\infty \alpha_n\frac{X^n}{n!}\right)=
\sum_{k=0}^\infty \sum_{n=0}^\infty P_{k,n}(\alpha_1,\dots,\alpha_n)
\frac{X^n}{n!}$$
and 
$$\log\left(1+\sum_{n=1}^\infty \alpha_n\frac{X^n}{n!}\right)=
\sum_{k=1}^\infty \sum_{n=0}^\infty (-1)^{k+1}(k-1)!
P_{k,n}(\alpha_1,\dots,\alpha_n)
\frac{X^n}{n!}$$
define the exponential function and the logarithm of 
an exponential generating series in $a\in\mathfrak m_{\mathcal E}$ 
respectively $1+a\in1+\mathfrak m_{\mathcal E}$ 
over an arbitrary  field $\mathbb K$.
These functions are one-to-one and mutually reciprocal.
\end{cor}

The following result shows that the functions $\exp!$
and $\log_!$ behave as expected 
under the derivation $\tau:\sum_{n=0}^\infty 
\alpha_n X^n\longmapsto \sum_{n=0}^\infty \alpha_{n+1}X^n$
of the shuffle-algebra.

\begin{prop} \label{propderivexplog} 
For all $A\in\mathfrak m=X\mathbb K[[X]]$ over an arbitrary field
$\mathbb K$ we have 
$$\tau\left({\exp}_!(A)\right)=({\exp}_!(A))\sh\tau(A)$$
and 
$$\tau\left({\log}_!(1+A)\right)=(1+A)^{\sh^{-1}}\sh\tau(A)$$
where $(1+A)^{\sh^{-1}}$ denotes the shuffle inverse of $(1+A)$.
\end{prop}

{\bf Proof} Proposition \ref{propshiftderiv} implies the formal 
identities
$$\tau\left(\sum_{n=0}^\infty \frac{A^{\sh^n}}{n!}\right)=
\sum_{n=0}^\infty n\frac{A^{\sh^{n-1}}}{n!}\sh\tau(A)=
\left(\sum_{n=0}^\infty
\frac{A^{\sh^n}}{n!}\right)\sh\tau(A)$$
for $A\in\mathfrak m$. By Proposition \ref{propcoeffintegr},
this identity holds over the ring $\mathbb Z$ and thus
over an arbitrary commutative field. This establishes
the formula for $\exp_!$.

For $\log_!$ we get similarly
$$\tau\left(-\sum_{n=1}^\infty \frac{(-A)^{\sh^n}}{n}\right)=
\sum_{n=1}^\infty n\frac{(-A)^{\sh^{n-1}}}{n}\sh\tau(A)$$
$$=
\left(\sum_{n=0}^\infty (-A)^{\sh^n}\right)\sh\tau(A)$$
which implies the result by Proposition \ref{propcoeffintegr}
and by the trivial identity 
$$(1+A)^{\sh^{-1}}=\sum_{n=0}^\infty (-A)^{\sh^n}$$
for the shuffle inverse $(1+A)^{\sh^{-1}}$ of $1+A\in 1+\mathfrak m$.
\hfill$\Box$

\section{The logarithm as a $p-$homogeneous form over 
$\overline{\mathbb F}_p[[x]]$}\label{sectlog}

Given a fixed prime number $p$,
Proposition \ref{propshgrFp} implies that there exists
polynomials $Q_{p,n}\in \mathbb N[\alpha_0,\dots,\alpha_n]$
for $n\geq 1$
such that 
$$\left(\sum_{n=0}^\infty \alpha_nX^n\right)^{\sh^p}=\alpha_0^p+
p\sum_{n=1}^\infty Q_{p,n}(\alpha_0,\dots,\alpha_n)X^n\ .$$
The polynomials $Q_{p,n}$ are homogeneous of degree $p$ with respect
to the variabels $\alpha_0,\dots,\alpha_n$ and 
we denote by 
$$\mu_p\left(\sum_{n=0}^\infty \alpha_nX^n\right)=\sum_{n=1}^\infty
Q_{p,n}(\alpha_0,\dots,\alpha_n)X^n$$
the $p-$homogeneous form defined by the ordinary generating
series of the polynomials $Q_{p,1},Q_{p,2},\dots$.

\begin{prop} \label{proploghom}
The restriction of $\mu_p$ to $1+\mathfrak m\subset
\overline{\mathbb F}_p[[X]]$ coincides with 
the function ${\log}_!$.
\end{prop}

{\bf Proof} We have 
$$\tau(\mu_p(1+A))=(1+A)^{\sh^{p-1}}\sh \tau(1+A)$$
for $A$ in $\mathfrak m$ where $\tau(\sum_{n=0}^\infty \alpha_nX^n)=
\sum_{n=0}^\infty \alpha_{n+1}X^n$ is the shift-operator 
of Proposition \ref{propshiftderiv}.
This identity defines the restriction of the
$p-$homogenous form $\mu_p$ to $1+\mathfrak m$. 
Proposition \ref{propderivexplog} and the identity 
$(1+A)^{{\sh}^{p-1}}\sh(1+A)=1$ show that the function ${\log}_!$
satisfies the same equation
$$\tau({\log}_!(1+A))=(1+A)^{{\sh}^{p-1}}\sh \tau(1+A)\ .$$

Since both series $\mu_p(1+A)$ and ${\log}_!(1+A)$ are
without constant term, the equality 
$\tau(\mu_p(1+A))=\tau({\log}_!(1+A))$ implies
the equality $\mu_p(1+A)={\log}_!(1+A)$.
\hfill$\Box$

\section{Proofs}\label{sectproofs}

\begin{prop}\label{propsra}
If $A$ in $X\overline{\mathbb F}_p[[X]]$
is rational (respectively algebraic)
then the formal power series ${\log}_!(1+A)$ is rational
(respectively algebraic).

More precisely, we have
$$\parallel {\log}_!(1+A)\parallel \leq 
1+{p+\parallel 1+A\parallel -1\choose p}\leq 1+\parallel 
1+A\parallel ^p$$
for $A$ rational in $\mathfrak m=X\overline{\mathbb F}_p[[X]]$,
respectively
$$\kappa( {\log}_!(1+A))\leq 
1+4\kappa(A){p+\kappa(1+A)-2\choose p-1}\leq
1+4(\kappa(1+A
))^p\ ,$$
for $A$ algebraic in $\mathfrak m$.
\end{prop}

\begin{prop}\label{propsira}
If $A$ in $X\overline{\mathbb F}_p[[X]]$
is rational (respectively algebraic)
then ${\exp}_!(A)$ is rational
(respectively algebraic).

More precisely, denoting by $q=p^e$ the cardinality of a finite field
$\mathbb F_q\subset \overline{\mathbb F}_p$ containing all coefficients
of $A$ we have
$$\parallel {\exp}_!(A)\parallel \leq p^{q^{\parallel A\parallel}}$$
for $A$ rational in $\mathfrak m$,
respectively
$$\kappa({\exp}_!(A)) \leq q^{\kappa(A)-1}p^{q^{\kappa(A)}}$$
for $A$ algebraic and non-zero in $\mathfrak m$.
\end{prop}

Theorems \ref{thmmainrat}, \ref{thmmainalg}, \ref{thmrateff}
and \ref{thmalgeff}
are now simple reformulations of  Propositions \ref{propsra}
and \ref{propsira}.

{\bf Proof of Proposition \ref{propsra}} 
The identity $(1+A)^{\sh^p}=1$ following from
Proposition \ref{propshgrFp} applied to the equality 
$$\tau({\log}_!(1+A))=(1+A)^{\sh^{-1}}\sh\tau(A)$$
of Proposition \ref{propderivexplog} establishes the equality
$$\tau({\log}_!(1+A))=(1+A)^{\sh^{p-1}}\sh\tau(A)$$
already encountered in the proof of Proposition \ref{proploghom}.
This shows 
$$\parallel \tau({\log}_!(1+A)\parallel \leq \parallel
1+A\parallel^{p-1}
\parallel \tau(A)\parallel \leq \parallel 1+A\parallel^p$$
and implies 
$$\parallel {\log}_!(1+A)\parallel \leq 1+\parallel 1+A
\parallel^p\ .$$
This proves the cruder inequality in the rational case.
The finer inequality follows from the fact that all $p$ factors
of $(1+A)^{\sh^{p-1}}\sh\tau(A)=\tau(\log_!(1+A))$ 
belong to a common vector space of dimension $\parallel 1+A\parallel$
which is closed for the shift map. The details are the same as for
Remark \ref{remrefinement}.

For algebraic $A$ we have similarly
$$\kappa( \tau({\log}_!(1+A)))\leq
(\kappa (1+A))^{p-1}\ \kappa (\tau(A))=
(\kappa (1+A))^{p-1}\ \kappa (\tau(1+A))\leq$$
$$\leq (\kappa(1+A))^{p-1}2\kappa(1+A)\leq
2(\kappa(1+A))^p$$
using assertion (i) of Proposition \ref{propKtau}. This shows
$$\kappa({\log}_!(1+A))\leq
1+2\kappa(\tau({\log}_!(1+A)))\leq
1+4(\kappa(1+A))^p$$
by assertion (ii) of Proposition \ref{propKtau}
and ends the proof for the cruder inequality.

The finer inequality follows from Proposition \ref{propKtau} 
combined with Remark \ref{remrefinement}.
\hfill$\Box$

Given a vector space
$\mathcal V\subset \mathbb K[[X]]$ containing $\mathbb K$, we denote by
$\Gamma(\mathcal V)$ the shuffle-subgroup generated by all elements of
$\mathcal V\cap (1+X\mathbb K[[X]])$.

\begin{lem} \label{lemVgrpe}
Every element of a vector space
$\mathcal V\subset \mathbb K[[X]]$ containing the field $\mathbb K$
of constants
can be written as a linear combination of elements in
$\Gamma(\mathcal V)$.
\end{lem}

{\bf Proof} We have the identity 
$$A=(1-\epsilon(A)+A)+(\epsilon(A)-1)$$
where $\epsilon(\sum_{n=0}^\infty \alpha_nX^n)=\alpha_0$ is the
augmentation map and where 
$(1-\epsilon(A)+A)$ and the constant 
$(\epsilon(A)-1)$ are both in $\mathbb K\Gamma(
\mathcal V)$ for $A\in\mathcal V$.\hfill$\Box$

{\bf Proof of Proposition \ref{propsira} for $A$ rational}
Corollary \ref{coralgfin}
shows that we can work over a finite subfield $\mathbb K=\mathbb F_q$ of
$\overline{\mathbb F}_p$ consisting of $q=p^e$ elements.

Given a rational series $A$ in $\mathfrak m=X\mathbb K[[X]]$, we 
denote by $\Gamma_A$ the shuffle-subgroup generated by all elements 
of the set 
$$\left\{\bigcup_{n=0}^\infty \left(\tau^n(A)+\mathbb K\right)\right\}
 \cap
\{ 1+X\mathbb K[[X]]\}.$$
This generating set of $\Gamma_A$ contains at most $q^{\parallel A
\parallel}$ elements. Proposition \ref{propshgrFp} implies thus
that $\Gamma_A$ is a finite group having at most $p^{q^{\parallel A
\parallel}}$ elements.
The subalgebra  
$\mathbb K[\Gamma_A]\subset \mathbb K[[X]]$
spanned by all elements of $\Gamma_A$ is thus of dimension
$\leq p^{q^{\parallel A\parallel}}$. The identity
$$\tau({\exp}_!(A))={\exp}_!(A)\sh\tau(A)$$
of Proposition \ref{propderivexplog}
and the fact that the derivation $\tau$ of
$\mathbb K[[X]]$ restricts to a derivation of the subalgebra
$\mathbb K[\Gamma_A]$ show the inclusion
$$\tau^n({\exp}_!(A))\in {\exp}_!(A)\sh\mathbb K[\Gamma_A]$$
for all $n\in\mathbb N$ by Lemma \ref{lemVgrpe}.
This ends the proof since the right-hand side is a 
$\mathbb K-$vector space of dimension at most $p^{q^{\parallel A\parallel}}$.
\hfill$\Box$

\begin{prop} We have for every prime number $p$ and for all natural integers
$j,k$ such that $j\geq 1$ the identity
$$\frac{(jk)!}{(j!)^kk!}\equiv \frac{(pjk)!}{((pj)!)^kk!}\pmod p\ .$$
\end{prop}

{\bf Proof} The number $(pjk)!/(((pj)!)^kk!)$ of the 
right-hand-side yields the cardinality of 
the set $\mathcal E$ of  all partitions of $\{1,\dots,pjk\}$
into $k$ subsets of $pj$ elements.
Consider the group $G$ generated by the $jk$ cycles of length $p$ of the
form $(i,i+jk,i+2jk,\dots,i+(p-1)jk)$ for $i=1,\dots,jk$.
The group $G$ has $p^{jk}$ elements and acts on the set of partitions 
by preserving their type defined as the multiset of cardinalities 
of all involved parts. In particular it acts by permutation on
the set $\mathcal E$. A partition $P\in \mathcal E$ 
is a fixpoint for $G$ if and only if every part of $P$ is
a union of $G-$orbits. Choosing a bijection between $\{1,\dots,jk\}$
and $G-$orbits of $\{1,\dots,pjk\}$,
fixpoints of $\mathcal E$ are 
in bijection with partitions of the set $\{1,\dots,jk\}$
into $k$ subsets of $j$ elements.
The number of fixpoints of the $G-$action on $\mathcal E$ equals thus
$\frac{(jk)!}{(j!)^kk!}$.
Since $G$ is a $p-$group, the cardinalities of all non-trivial 
$G-$orbits of $\mathcal E$ are strictly positive powers of $p$. This ends
the proof.
\hfill$\Box$

\begin{cor} ${\exp}_!$ and ${\log}_!$ commute with
the \lq\lq Frobenius substitution\rq\rq 
$$\varphi(\sum_{n=0}^\infty \alpha_nX^n)=\sum_{n=0}^\infty \alpha_nX^{pn}$$
for series in $X\overline{\mathbb F}_p[[X]]$, respectively 
in $1+X\overline{\mathbb F}_p[[X]]$.
\end{cor}

This implies 
$$({\exp}_!(A))_{0,f}={\exp}_!(A_{0,f})$$
where $C_{k,f}=\sum_{n=0}^\infty 
\gamma_{k+np^f}X^n$ for $C=\sum_{n=0}^\infty \gamma_nX^n$.

\begin{lem}\label{lem01}
We have 
$$(B\sh C)_{0,1}=B_{0,1}\sh C_{0,1}$$
\end{lem}

{\bf Proof} Follows from the identity 
$${pn\choose k}\equiv 0\pmod p$$ if 
$k\not\equiv 0\pmod p$.\hfill$\Box$

{\bf Proof of Proposition \ref{propsira} for $A$ algebraic}
We work again over a finite subfield $\mathbb K=\mathbb F_q\subset 
\overline{\mathbb F}_p$ containing all coefficients of $A$.

Let $\Gamma_A$ denote the shuffle-subgroup generated by 
all elements in 
$$(\mathcal K(A)+\mathbb K)\cap(1+X\mathbb K[[X]])$$
where $\mathcal K(A)=\mathbb KA+\sum_{k,f}\mathbb K A_{k,f}$
denotes the $p-$kernel of $A$.
We denote by $\mathbb K[\Gamma_A]\subset (\mathbb K[[X]],\sh)$
the shuffle-subalgebra of dimension at most $p^{q^{\kappa(A)}}$
spanned by all elements of the group $\Gamma_A\subset (1+X\mathbb K[[X]],
\sh)$.

Using the convention $A_{0,0}=A$, we have for $B\in\mathbb K[\Gamma(A)]$
and for $k$ such that $0\leq k<p$
$$\left({\exp}_!(A_{0,f})\sh B\right)_{k,1}=\left(\tau^k\left(
{\exp}_!(A_{0,f})\sh B\right)\right)_{0,1}=$$
$$=\left(\sum_{j=0}^k{k\choose j}\tau^j({\exp}_!(A_{0,f}))\sh
\tau^{k-j}(B)\right)_{0,1}=$$
$$=\sum_{j=0}^k{k\choose j}
\left(\tau^j({\exp}_!(A_{0,f}))\right)_{0,1}\sh B_{k-j,1}$$
where the last equality is due to Lemma \ref{lem01} (and to the
equality $(\tau^k(C))_{0,1}=C_{k,1}$ for $0\leq k<p$).

Iteration of the identity 
$\tau({\exp}_!(A_{0,f}))={\exp}_!(A_{0,f})
\sh \tau(A_{0,f})$ given by Proposition \ref{propderivexplog}
shows that $\tau^j({\exp}_!(A_{0,f}))$ is of the form
${\exp}_!(A_{0,f})\sh F$ where $F$ is a linear combination
of shuffle-products involving at most $j$ factors of the set 
$\{\tau(A_{0,f}),\tau^2(A_{0,f}),\dots,\tau^j(A_{0,f})\}$.
Applying Lemma \ref{lem01} we get
$$(\tau^j({\exp}_!(A_{0,f})))_{0,1}=(
{\exp}_{0,f+1}(A))\sh F_{0,1}\ .$$
An iterated application of Lemma \ref{lem01} shows now that
$F_{0,1}$ is a linear combination of shuffle-products involving
at most $j$ factors in $\{A_{1,f+1},\dots,A_{j,f+1}\}$. We have thus
$F_{0,1}\in \mathbb K[\Gamma_A]$ by Lemma \ref{lemVgrpe} and we
get the inclusion
$$({\exp}_!(A_{0,f})\sh\mathbb K[\Gamma_A])_{k,1}\subset
{\exp}_!(A_{0,f+1})\sh \mathbb K[\Gamma_A]$$
for all $f\in\mathbb N$ and for all $k\in\{0,\dots,p-1\}$.

Setting
$$E_A=\{{\exp}_!(B)\ \vert\ B\in\mathcal K(A)\cap X\mathbb K[[X]]\}$$
we have the inclusion
$$\mathcal K({\exp}_!(A))\subset 
E_A\sh \mathbb K[\Gamma_A]
\subset\mathbb K[E_A]\sh \mathbb K[\Gamma_A]$$
where $\mathcal K(\exp_!(A))$ denotes the $p-$kernel 
of $\exp_!(A)$.
This implies 
$$\kappa({\exp}_!(A))\leq
\dim(\mathbb K[E_A)\dim(\mathbb K[\Gamma_A])
\ .$$

We suppose now $A$ non-zero.
The vector space $\mathcal K(A)\cap X\mathbb K[[X]]$ is thus of 
codimension $1$ in $\mathcal K(A)$. The image $E_A$ of 
$\mathcal K(A)\cap X\mathbb K[[X]]$ under the group-isomorphism
$\exp_!:(X\mathbb K[[X]],+)\longmapsto (1+X\mathbb K[[X]],\sh)$ 
is hence a subgroup of cardinality $q^{\kappa(A)-1}$ in 
$(1+X\mathbb K[[X]],\sh)$.
We have thus
$$\kappa({\exp}_!(A))\leq
\dim(\mathbb K[E_A])\dim(\mathbb K[\Gamma_A])\leq
q^{\kappa(A)-1}p^{q^{\kappa(A)}}$$
which ends the proof.\hfill$\Box$.

\section{Power series in free non-commuting variables}
\label{sectionpows}

This and the next section recall a 
few basic and well-known facts concerning
(rational) power series in free non-commuting variables, 
see for instance \cite{St2}, \cite{BR} or a similar book
on the subject. We use however sometimes a different terminology,
motivated by \cite{B}.

We denote by $\mathcal X^*$ the free monoid on a finite set $\mathcal X=
\{X_1,\dots,X_k\}$. We write $1$ for the identity element 
and we use a boldface capital 
$\mathbf X$ for a non-commutative monomial
$\mathbf X=X_{i_1}X_{i_2}\cdots X_{i_l}\in \mathcal X^*$. 
We denote by
$$A=\sum_{\mathbf X\in\mathcal X^*}(A,\mathbf X)\mathbf X\in
\mathbb K\dlangle X_1,\cdots ,X_k\drangle$$
a non-commutative formal power series where
$$\mathcal X^*\ni\mathbf X\longmapsto (A,\mathbf X)\in \mathbb K$$
stands for the coefficient function.

We denote by 
$\mathfrak m\subset \mathbb K\dlangle X_1,\dots,X_k\drangle $
the maximal ideal consisting of formal power series without
constant coefficient and by 
$\mathbb K^*+\mathfrak m=\mathbb K\dlangle X_1,\dots,X_k\drangle
\setminus\mathfrak m$ the unit-group of the
algebra $\mathbb K\dlangle X_1,\dots,X_k\drangle$
consisting of all (multiplicatively) invertible elements in
$\mathbb K\dlangle X_1,\dots,X_k\drangle $. 
The unit group is
isomorphic to the direct product $\mathbb K^*\times (1+\mathfrak m)$
where $\mathbb K^*$ is the central subgroup consisting of
non-zero constants and where $1+\mathfrak m$ 
denotes the multiplicative subgroup given by the affine subspace
formed by power series with
constant coefficient $1$. We have $(1-A)^{-1}=1+\sum_{n=1}^\infty
A^n$ for the multiplicative inverse $(1-A)^{-1}$
of an element $1-A\in 1+\mathfrak m$. 

\subsection{The shuffle algebra}\label{sectshufflealg}

The {\it shuffle-product} $\mathbf X\sh \mathbf X'$ of two 
non-commutative monomials
$\mathbf X,\mathbf X'\in \mathcal X^*$ of degrees 
$a=\mathop{deg}(\mathbf X)$ and $b=\mathop{deg}(\mathbf X')$ 
(for the obvious
grading given by $\mathop{deg}(X_1)=\dots=\mathop{deg}(X_k)=1$)
is the sum
of all ${a+b\choose a}$ monomials of degree $a+b$ obtained by 
``shuffling'' in all possible ways the linear factors (elements of
$\mathcal X$) involved in $\mathbf X$ 
with the linear factors of $\mathbf X'$. A monomial involved in
$\mathbf X\sh \mathbf X'$ can be thought of as a monomial 
of degree $a+b$ whose
linear factors are coloured by two colours with $\mathbf X$ corresponding 
to the product of all linear factors of the first colour and
$\mathbf X'$ corresponding to the product of the remaining 
linear factors. The shuffle product $\mathbf X\sh \mathbf X'$ 
can also be recursively defined by $\mathbf X\sh 1=1\sh \mathbf X=
\mathbf X$ and 
$$(\mathbf X X_s)\sh(\mathbf X' X_t)=(\mathbf X\sh(\mathbf X' X_t))
X_s+((\mathbf X X_s)\sh
\mathbf X')X_t$$ where
$X_s,X_t\in\mathcal X=\{X_1,\dots,X_k\}$ are monomials of degree $1$.

Extending 
the shuffle-product in the obvious way to formal power series endows the
vector space $\mathbb K\dlangle X_1,\dots,X_k\drangle $ with an associative and
commutative algebra structure called the
{\it shuffle-algebra}. In the case of one variable $X=X_1$ we recover the
definition of Section \ref{sectshufflpr}.

The group $\hbox{GL}_k(\mathbb K)$ acts 
on the vector space $\mathbb K\dlangle X_1,\dots,
X_k\drangle $ by linear substitutions. 
This action induces an automorphism of the multiplicative
(non-commutative) algebra-structure or of
the (commutative) shuffle algebra-structure underlying 
the vector space $\mathbb K\dlangle X_1,\dots,X_k\drangle$. 

Substitution of all variables $X_j$ of formal power series in 
$\mathbb K\dlangle X_1,\dots,X_k\drangle $ by $X$ (or more 
generally by arbitrary not necessarily equal
formal power series without constant term)
yields a homomorphism of (shuffle-)algebras into the commutative 
(shuffle-)algebra $\mathbb K[[X]]$.

The commutative unit group (set of invertible elements for the
shuffle-product) of the shuffle algebra, given by the set
$\mathbb K^*+\mathfrak m$, is isomorphic to the direct product
$\mathbb K^*\times (1+\mathfrak m)$ where $1+\mathfrak m$ is endowed with the 
shuffle product. The inverse of an element
$1-A\in(1+\mathfrak m,\sh)$ is given by $\sum_{n=0}^\infty A^{\sh^n}=
1+A+A\sh A+A\sh A\sh A+\dots$. 

The following result generalises Proposition \ref{propshgrFp}:

\begin{prop} \label{propshgrFpgen} Over a field of positive characteristic $p$,
the subgroup $1+\mathfrak m$ of the shuffle-group is an
infinite-dimensional $\mathbb F_p-$vector space.
\end{prop}

{\bf Proof} Contributions to a $p-$fold shuffle product
$A_1\sh A_2\sh\cdots \sh A_p$ are given by 
monomials with linear factors coloured by 
$p$ colours $\{1,\dots,p\}$
keeping track of their ``origin'' with coefficients given
by the product of the corresponding ``monochromatic'' coefficients
in $A_1,\dots,A_p$. A permutation of 
the colours $\{1,\dots,p\}$ (and in particular, 
a cyclic permutation of all 
colours) leaves such a contribution invariant if
$A_1=\dots=A_p$. Coefficients of strictly
positive degree in $A^{\sh^p}$ are thus zero in characteristic $p$.
\hfill$\Box$

As in the one variable case, one can prove that 
$$\frac{1}{k!}A^{{\sh}^k}$$
is defined over an arbitrary field $\mathbb K$ for $A\in\mathfrak m$.
Indeed, monomials contributing to $A^{{\sh}^k}$ can be considered as 
colored by $k$ colours and the $k!$ possible colour-permutations
yield identical contributions.

For $A\in\mathfrak m$, we denote by 
$${\exp}_!(A)=\sum_{n=0}^\infty \frac{1}{n!}A^{{\sh}^n}$$
the resulting exponential map from the Lie algebra $\mathfrak m$
into the infinite-dimensional commutatif Lie group $(1+\mathfrak m,\sh)$.
As expected, its reciprocal function is defined by
$${\log}_!(1+A)=\sum_{n=1}^\infty \frac{(-1)^{n+1}}{n}A^{{\sh}^n}\ .$$
In the case of a field $\mathbb K$ of positive characteristic $p$
the function ${\log}_!$ is again given by the 
restriction to $1+\mathfrak m$ of a $p-$homogeneous form 
$\mu_p$. 

The form $\mu_p$ has all its coefficients in $\mathbb N$ and
is again defined by the equality 
$$A^{{\sh}^p}=(A,1)^p+p\mu_p(A)$$
over $\mathbb Z$. It can thus be defined over an arbitrary field.

\section{Rational series}

A formal power series $A$ is {\it rational} if it belongs to
the smallest 
subalgebra in $\mathbb K\dlangle X_1,\dots,X_k\drangle $
which contains the free
associative algebra $\mathbb K\langle
X_1,\dots,X_k\rangle$ of non-commutative polynomials
and intersects the multiplicative unit group of
$\mathbb K\dlangle X_1,\dots,X_k\drangle $ in a subgroup.

Given a monomial $\mathbf T\in\mathcal X^*$, we
denote by $$\rho(\mathbf T):\mathbb K \dlangle X_1,\dots,X_k\drangle 
\longrightarrow \mathbb K \dlangle X_1,\dots,X_k\drangle $$ 
the linear application defined by 
$$\rho(\mathbf T)A=\sum_{\mathbf X\in\mathcal X^*}(A,\mathbf X\mathbf
T)\mathbf X$$
for $A=\sum_{\mathbf X\in\mathcal X^*}(A,\mathbf X)\mathbf X$ in
$\mathbb K\dlangle X_1,\dots,X_k\drangle$.
The identity $\rho(\mathbf T)(\rho(\mathbf T')A)=\rho(
\mathbf T\mathbf T')A$ shows that we have a representation
$$\rho:\mathcal X^*\longrightarrow\hbox{End}(\mathbb K\dlangle \mathcal X
\drangle)$$
of the free monoid $\mathcal X^*$ on $\mathcal X$. 
The {\it recursive closure} $\overline A$ of a power series $A$
is the vector space
spanned by its orbit $\rho(\mathcal X^*)A$ under $\rho(\mathcal X^*)$.
We call the dimension $\dim(\overline A)$ of $\overline A$ the
{\it complexity} of $A$.

We call a subspace
$\mathcal A\subset \mathbb K\dlangle X_1,\dots,X_k\drangle $ {\it recursively
  closed} if it contains the recursive closure of all its elements.

Rational series coincide with series of finite complexity by a 
Theorem of Sch\"utzenberger (cf. \cite{BR}, Theorem 1 of page 22).

\begin{rem} \label{remdimratfrac}
In the case of one variable, the complexity 
$\mathop{dim}(\overline A)$ of a reduced non-zero 
rational fraction $A=\frac{f}{g}$ with
$f\in\mathbb K[X]$ and $g\in1+X\mathbb K[X]$ equals
$\mathop{dim}(\overline
A)=\mathop{max}(1+\mathop{deg}(f),\mathop{deg}(g))$.
\end{rem}

\begin{rem}
The (generalised) {\it Hankel matrix} $H=H(A)$ of
$$A=\sum_{\mathbf X\in\mathcal X^*}(A,\mathbf X)\mathbf X\in
\mathbb K\dlangle X_1,\dots,X_k\drangle $$ 
is the infinite matrix with rows and columns 
indexed by the free monoid $\mathcal X^*$ of monomials and entries
$H_{\mathbf X\mathbf X'}=(A,\mathbf X\mathbf X')$. 
The rank $\mathop{rank}(H)$ is given by the
complexity $\dim(\overline A)$ of $A$ and 
$\overline A$ corresponds to the row-span of $H$.
\end{rem}

Given subspaces $\mathcal A,\mathcal B$ of $\mathbb K\dlangle\mathcal X
\drangle$, we denote by $\mathcal A\sh \mathcal B$ the vector space
spanned by all products $A\sh B$ with $A\in\mathcal A$ and 
$B\in\mathcal B$.

\begin{prop} \label{propratclosed}
We have the inclusion
$$\overline{A\sh B}\subset \overline A\sh \overline B$$
for the shuffle product $A\sh B$ of $A,B\in\mathbb K\dlangle X_1,\dots,
X_k\drangle $.
\end{prop}

The following result is Proposition 4 of \cite{Fl}:

\begin{cor} \label{corratclosed}
We have
$$\mathop{dim}(\overline{A\sh B})\leq 
\mathop{dim}(\overline A)\ \mathop{dim}(\overline B)$$
for the shuffle product $A\sh B$ of $A,B\in\mathbb K\dlangle X_1,\dots,
X_k\drangle $.

In particular, shuffle products of rational elements in $\mathbb
K\dlangle X_1,\dots,X_k\drangle $ are rational.
\end{cor}

{\bf Proof of Proposition \ref{propratclosed}} 
For $Y\in\overline A,Z\in 
\overline B$ and $X$ in $\{X_1,\dots,X_k\}$, 
the recursive definition of the shuffle product given in Section
\ref{sectshufflealg} shows
$$\rho(X)(Y\sh Z)=(\rho(X)Y)\sh Z+ Y\sh (\rho(X)Z)\ .$$
We have thus the inclusions
$$\rho(X)(Y\sh Z)\in \overline A\sh Z+Y\sh \overline B
\subset \overline A\sh\overline B$$
which show that the vector space
$\overline A\sh \overline B$ is recursively closed.
Proposition \ref{propratclosed} follows now from the inclusion
$A\sh B\in \overline A\sh \overline B$.
\hfill$\Box$

\begin{rem} Similar arguments show that the set of 
rational elements in $\mathbb K\dlangle X_1,\dots,X_k\drangle$
is also closed under the ordinary product
(and multiplicative inversion of invertible series),
Hadamard product and composition (where
one considers $A\circ(B_1,\dots,B_k)$ with $A\in\mathbb
K\dlangle X_1,\dots,X_k\drangle $ and $B_1,\dots,B_k\in\mathfrak m\subset 
\mathbb K\dlangle X_1,\dots, X_k\drangle $).
\end{rem}

\begin{rem} The shuffle inverse of a rational element in 
$\mathbb K^*+\mathfrak m$ is in general not rational in characteristic
$0$. An exception is given by geometric progressions
$\frac{1}{1-\sum_{j=1}^k \lambda_j X_j}=\sum_{n=0}^\infty 
\left(\sum_{j=1}^k\lambda_j X_j\right)^n$ since we have
$$\frac{1}{1-\sum_{j=1}^k \lambda_j X_j}\sh
\frac{1}{1-\sum_{j=1}^k \mu_j X_j}=
\frac{1}{1-\sum_{j=1}^k (\lambda_j+\mu_j) X_j}$$
corresponding to $e^{\lambda X}e^{\mu
  X}=e^{(\lambda+\mu)X}$ 
in the one-variable case.

There are no other such 
elements in $1+\mathfrak m\subset \mathbb K[[X]]$, see 
Remark \ref{remshfflegrpnoalgelts}.
I ignore if the maximal rational shuffle subgroup
of $1+\mathfrak m\subset \mathbb C\dlangle X_1,
\dots,X_k\drangle $ (defined as the set
of all rational elements in $1+\mathfrak m$ with rational inverse
for the shuffle product) contains 
other elements if $k\geq 2$.
\end{rem}

\begin{rem} Any
finite set of rational elements in $\mathbb K\dlangle X_1,
\dots,X_k\drangle$ over a field $\mathbb K$ of positive 
characteristic is included in
a unique minimal finite-dimensional recursively closed subspace
of $\mathbb K\dlangle X_1,\dots,X_k\drangle $ which 
intersects the shuffle group
$(\mathbb K^*+\mathfrak m,\sh)$ in a subgroup.  
\end{rem}

\section{Main result for generating series in non-commuting variables
}\label{sectnonlin}

The following statement is our main result in a non-commutative
frame-work.

\begin{thm} Let $\mathbb K$ be a subfield of $\overline{\mathbb F}_p$.
Given a non-commutative formal power 
series $A\in\mathfrak m\subset
\mathbb K\dlangle \mathcal X\drangle$, the following two assertions
are equivalent:

\ \ (i) $A$ is rational.

\ \ (ii) ${\exp}_!(A)$ is rational.

More precisely, we have for a rational series $A$ in 
$\mathfrak m$ the inequalities
$$\dim\left(\overline{{\log}_!(1+A)}\right)\leq
1+\left(\dim(\overline{1+A})\right)^p$$
and 
$$\dim\left(\overline{{\exp}_!(A)}\right)\leq
p^{q^{\dim(\overline{A})}}$$
where $q=p^e$ is the cardinality of a finite field 
$\mathbb F_q$ containing all coefficients of $A$.
\end{thm}

{\bf Proof} The identity 
$${\log}_!(1+A)=\sum_{X\in\mathcal X}\left((1+A)^{{\sh}^{p-1}}
\sh \rho(X)A\right)X$$
and Corollary \ref{corratclosed} show
$$\dim\left(\overline{{\log}_!(1+A)}\right)\leq
1+\left(\dim(\overline{1+A})\right)^p\ .$$

For the opposite direction we denote by $\mathbb K=\mathbb F_q$
a finite subfield of $\overline{\mathbb F}_p$ containing all 
coefficients of $A$. We have
$$\overline{{\exp}_!(A)}\subset {\exp}_!(A)\sh
\mathbb K[\Gamma(A)]$$
where $\mathbb K[\Gamma(A)]$ is the shuffle subalgebra of dimension
$\leq p^{q^{\dim(\overline A)}}$ spanned by all elements of the group
$\Gamma$ generated by 
all elements of the form
$$(\overline A+\mathbb K)\cap(1+\mathfrak m)\ .$$
This implies the inequality
$$\dim\left(\overline{{\exp}_!(A)}\right)\leq p^{q^{
\dim(\overline A)}}$$
which ends the proof.\hfill$\Box$

{\bf Acknowledgements} I thank J-P. Allouche, 
M. Brion, A. Pantchichkine, T. Rivoal, B. Venkov and an anonymous referee
for their interest and helpful remarks.

\noindent Roland BACHER

\noindent INSTITUT FOURIER

\noindent Laboratoire de Math\'ematiques

\noindent UMR 5582 (UJF-CNRS)

\noindent BP 74

\noindent 38402 St Martin d'H\`eres Cedex (France)
\medskip

\noindent e-mail: Roland.Bacher@ujf-grenoble.fr

\end{document}